\documentclass[a4paper]{article}
\usepackage{amsmath,amssymb,amsthm,eucal,graphics}
\usepackage[cp1251]{inputenc}
\usepackage[russian]{babel}
\usepackage{graphicx}
\usepackage{amsfonts,amssymb,eucal}
\usepackage{amsbsy}
\usepackage{latexsym}
\usepackage{textcase}
\usepackage{cite}
\usepackage{enumerate}
\usepackage{cmap}
\usepackage{comment}

\usepackage[pdftex,unicode,colorlinks,linkcolor=blue,
citecolor=red,bookmarksopen,pdfhighlight=/N]{hyperref}

\newtheorem{theorem}{Теорема}

\newtheorem{lemma}{Лемма}

\newtheorem*{theoBMm}{Теорема Бёрлинга\,--\,Мальявена о мультипликаторе}

\theoremstyle{definition}

\newtheorem{rem}{Замечание}
\newtheorem{example}{Пример}

\renewcommand{\leq}{\leqslant} 
\renewcommand{\geq}{\geqslant}
\newcommand{\RR}{\mathbb{R}} 
\newcommand{\CC}{\mathbb{C}} 
\newcommand{\NN}{\mathbb{N}} 

\newcommand{\DD}{\mathbb{D}}

\DeclareMathOperator{\bal}{bal}

\DeclareMathOperator{\type}{type} 
\DeclareMathOperator{\strip}{str} 

 \DeclareMathOperator{\up}{up}
\DeclareMathOperator{\dd}{d}

\renewcommand{\Re}{\operatorname{Re}}
\renewcommand{\Im}{\operatorname{Im}}

\begin{document}

\Russian
\sloppy
\rm

\begin{flushright}
УДК 517.53+517.574
\end{flushright}

\begin{large}
\begin{center}

\textbf{Субгармоническое дополнение к теореме  Бёрлинга\,--\,Мальявена о мультипликаторе}

\

{Б.\,Н. Хабибуллин, Е.\,Г. Кудашева}\footnote{Работа 
%%%поддержана грантом Российского научного фонда (проект № 22-21-00026) и 
выполнена в рамках реализации программы развития Научно-образовательного математического центра Приволжского федерального округа (соглашение №  075-02-2022-882).}
\end{center}
\end{large}

\medskip
\begin{abstract}
Теорема Бёрлинга\,--\,Мальявена о мультипликаторе и различные её версии 
дают несколько вариантов условий на  функцию $f$ на вещественной оси $\mathbb R$, при которых 
эту функцию можно умножить на ограниченную на $\mathbb R$ целую функцию $h$ сколь угодно малого зкспоненциального типа   $>0$  так, что  произведение $fh$ ограничено на $\mathbb R$.  
Мы рассматриваем новую версию для функций $f=\exp(u-M)$,  где $u$ и $M$ --- пара субгармонических функций конечного типа с конечными логарифмическими интегралами по $\RR$.

\begin{comment}
Subharmonic addition to the Beurling\,--\,Malliavin  Theorem on the multiplier

We prove the a version of the Beurling\,--\,Malliavin Theorem on multiplier. This version is formulated here in a simplified form.
 Let $u\not\equiv -\infty$ and $M\not\equiv -\infty$ be a pair of subharmonic   functions on the complex plane $\mathbb C$ with positive parts $u^+:=\sup\{u,0\}$ and $M^+$  
such that
$$
\operatorname{type}[u]:=\limsup_{z\to \infty} \frac{u^+(z)}{|z|}<+\infty, \;
\operatorname{type}[M]<+\infty, \;    
\int_{-\infty}^{+\infty}\frac{u^+(x)+M^+(x)}{1+x^2}\operatorname{d}x<+\infty.
$$
If  $\operatorname{type}[u]<a<+\infty$, $0<b<+\infty$, and $\operatorname{type}[M]<c<+\infty$, 
 then there are an entire function $h\not\equiv 0$ with  $\operatorname{type}[\log|f|]<c$ and 
a subset $iY$ in the imaginary axis of linear Lebesgue measure  $<b$ such that the function $h$ is ouncec 
$u(z)-M(z)+\log|h(z)|\leq a|\Im z|$
on each straight line parallel to the real axis and not intersecting $iY$. 

\end{comment}

\medskip
\noindent \textit{Ключевые слова}: 
целая функция экспоненциального типа, субгармоническая функция конечного типа,   
мультипликатор, класс Картрайт 
\end{abstract}

\section{Введение и  основной результат}

Одноточечные множества $\{x\}$ часто записываем без фигурных скобок, т.е. просто как $x$. Так,   $\NN_0:=0\cup \NN=\{0,1, \dots\}$ для  множества $\mathbb N:=\{1,2, \dots\}$ {\it натуральных чисел.\/}  
Через  $\CC$ и $\RR$ обозначаем соответственно \textit{комплексную плоскость\/} и 
\textit{действительную прямую,\/}     часто рассматриваемую ниже как \textit{вещественную ось\/}
$\RR\subset \CC$, с  {\it расширением\/}    
двумя <<бесконечными значениями>> $-\infty:=\inf \RR\notin \RR$, $+\infty:=\sup \RR\notin \RR$, неравенствами $-\infty\leq x\leq +\infty$ для любого $x\in \overline \RR$
 и  порядковой топологией с базой открытых множеств из открытых интервалов 
$(a,b):=\bigl\{x\in \overline \RR\bigm| a<x<b \bigr\}$ при  $a<b$, а также $[-\infty,b):=(-\infty,b)\cup -\infty$
 и $(a,+\infty]:=(a,+\infty)\cup +\infty$.  
Символом  $0$, кроме нуля,   могут обозначаться {\it нулевые\/} функции,  меры  и пр., а через $-\infty$
или $+\infty$ и функции, тождественно равные $-\infty$ или $+\infty$. 
Если  $a\in \overline{\RR}$ или $a\colon X\to \overline{\RR}$ --- {\it расширенная числовая\/} функция,  то 
$a^+:=\sup\{a,0\}$  --- \textit{положительная часть\/} соответственно числа или функции $a$, а для подмножества 
$A\subset \overline{\RR}$ полагаем $A^+:=\bigl\{a^+\bigm| a\in A\bigr\}$.
К примеру, $\RR^+=[0,+\infty)$ --- \textit{положительная полуось} с \textit{расширением\/}  
$\overline{\RR}^+=\RR^+\cup +\infty$.

 Для  голоморфной на $\CC$, или \textit{целой,\/} функции $f$
\begin{equation}\label{efet}
\type_f:=\limsup_{z\to \infty}\frac{\ln^+|f(z)|}{|z|}\in \overline{\RR}^+
\end{equation}
--- величина её верхнего типа  при порядке $1$, или далее просто \textit{тип\/} целой функции $f$. 
Если при этом $\type_f\in \RR^+$,  то функция $f$ называется {\it целой функцией   экспоненциального типа\/} 
\cite{Boas}, \cite{Levin96}, \cite{RC},  \cite{Khsur}, хотя так же  широко распространён 
и термин <<целая функция конечной степени>> \cite{Levin56}, \cite{HJ94}. 
Интеграл от функции $v\colon \RR\to \overline \RR$
\begin{equation}\label{Jm0}
J[v]:=\frac{1}{\pi}\int_{-\infty}^{+\infty}\frac{v(x)}{1+x^2}\dd x
\end{equation}
часто называют \textit{логарифмическим интегралом} \cite{Koosis88}, \cite{Koosis92}, \cite{Koosis96}
функции $v$. 
 
\begin{theoBMm}[{\cite{BM62}, \cite{Mal79}, 
\cite{Koosis88}, \cite{Koosis92}, \cite{Koosis96}, \cite{HJ94}, \cite{MasNazHav05}, \cite{Khsur}}]  Пусть функция $u\colon \RR\to \RR^+$ равномерно непрерывна на $\RR$
или же совпадает с сужением на $\RR$  функции $\ln |F|$, где $F\neq 0$ --- целая функция экспоненциального типа.  
Если логарифмический интеграл $J[u^+]$  конечен, 
то для любого числа $c> 0$ существует такая ограниченная на $\RR$ целая функция $h$ экспоненциального типа $\type_h\leq c$, что имеет место неравенство 
\begin{equation}\label{uh}
u(x)+\ln \bigl|h(x)\bigr|\leq 0 \quad\text{при всех $x\in \RR$}.  
\end{equation}
\end{theoBMm} 
В случае $u=\ln|f|$ неравенство \eqref{uh} чаще и более традиционно записывают как 
$|fh|\leq 1$ на $\RR$, что, по-видимому, и обусловило название этого замечательного результата как теоремы о мультипликаторе, т.е. множителе $h$, <<гасящем рост>> $f$ до ограниченности произведения $fh$ на $\RR$. 
Мы дополним эту теорему некоторой новой версией,  в которой  основную роль играют  субгармонические функции. 

Для расширенной числовой функции $u$ на $\CC$ величина   
\begin{equation} \label{typevf}
\type[u]:=\limsup_{z\to \infty} \frac{u^+(z)}{|z|}\in \overline \RR
 \end{equation} 
(верхний) {\it --- тип\/} (роста) функции $u$ {\it при порядке\/} $1$ (около $+\infty$) \cite{Boas}, \cite{Levin56}, \cite{Levin96}, \cite{Kiselman}, \cite[2.1]{KhaShm19}, или просто \textit{тип\/} функции $u$ 
без упоминания порядка $1$ далее.
Функции $u$ конечного типа, если  $\type[u]\in \RR^+$. К примеру, для  \eqref{efet} имеем $\type_f=\type\bigl[\ln |f|\bigr]$. 

 Далее $D_z(r):=\bigl\{w \in \CC \bigm| |w-z|<r\bigr\}$ и  $\overline{D}_z(r):=\bigl\{w \in \CC_{\infty} \bigm| |w-z|\leq r\bigr\}$,
а также $\partial \overline{D}_z(r):=\overline{D}_z(r)\setminus {D}_z(r)$
---  соответственно {\it открытый} и   {\it замкнутый круги,\/} а также  {\it окружность  радиуса\/ $r\in \overline \RR^+$
 с центром\/} $z\in \CC$, а   $\DD:=D_0(1)$ и  $\overline \DD:=\overline D_0(1)$, а также  
$\partial \overline \DD:=\partial \overline{D}_0(1)$ --- соответственно {\it открытый\/} и {\it замкнутый 
единичные круги,\/} а также {\it единичная окружность} в  $\CC$. 

Через $\CC^{\up}:=\bigl\{z\in \CC \bigm| \Im z>0\bigr\}$ и $\CC^{\overline  \up}:=\CC^{\up}\cup \RR$ обозначаем 
{\it верхние\/} соответственно {\it открытую} и   {\it замкнутую полуплоскости,\/} 
а $-\CC^{\up}$ и $-\CC^{\overline  \up}$ --- это  соответственно {\it нижние   открытая} и   {\it замкнутая полуплоскости\/} в  $\CC$. 

{\it Сужение\/} функции или меры $m$ на  $S\subset \CC$ обозначаем как   $m{\lfloor}_S$.

Следуя \cite[определение 3]{Kha22}, для $d\in  \RR^+$, полунерерывной снизу функции  $r\colon \CC\to   \overline \RR^+\setminus 0$ и {\it гамма-функции\/} $\Gamma$ внешнюю   меру 
 \begin{equation}\label{mr}
{\mathfrak m}_d^r\colon S\underset{S\subset \CC}{\longmapsto}  
\inf \Biggl\{\sum_{k} \dfrac{\pi^{d/2}}{\Gamma (1+d/2)}r_k^d\biggm| S\subset \bigcup_{k} 
\overline D_{z_k}(r_k), \, z_k\in \CC, \, r_k \leq  r(z_k)\Biggr\}
\end{equation}
 называем {\it $d$-мерным обхватом Хаусдорфа переменного радиуса обхвата\/ $r$.\/} При этом  
через {\it постоянные функции\/} $r>0$ определяется {\it $d$-мерная    мера Хаусдорфа\/}
\begin{equation*}
{\mathfrak m}_d\colon S\underset{S\subset \CC}{\longmapsto}  \lim_{0<r\to 0} {\mathfrak m}_d^r(S)
\underset{r>0}{\geq} {\mathfrak m}_d^r(S)\geq {\mathfrak m}_d^\infty(S),
\end{equation*}
являющаяся регулярной  мерой Бореля, и  
${\mathfrak m}_d\geq {\mathfrak m}_d^{r}\geq {\mathfrak m}_d^{t}\geq {\mathfrak m}_d^{\infty}$
для любых пар  функций $r\leq t$.
В частности,   ${\mathfrak m}_2$ ---  это {\it плоская мера  Лебега на\/} $\CC$,  а для любой  липшицевой кривой $L$ в $\CC$  сужение   ${\mathfrak m}_1{\lfloor}_L$ --- это мера длины дуги на липшицевой 
кривой $L$ \cite[3.3.4A]{EG}. Таким образом,   ${\mathfrak m}_1{\lfloor}_{\RR}$ --- это обычная \textit{линейная  мера  Лебега на\/} $\RR$.
Кроме того,  $0$-мерная   мера Хаусдорфа ${\mathfrak m}_0$ множества -- это число элементов в нём, а также 
${\mathfrak m}_d={\mathfrak m}_d^r=0$ при любом    $d>2$.

Как и  в \cite[предисловие]{EG}, расширенная числовая функция {\it интегрируема\/} по мере  Бореля  $\mu$, 
или \textit{$\mu$-интегрируема,\/} если интеграл от неё по этой мере  корректно определён значением из $\overline \RR$, и  {\it суммируема\/} по   $\mu$, или  {\it $\mu$-суммируема,\/}  если этот   интеграл конечен, т.е. принимает значения из $\RR$.  

Рассмотрим   функцию  $r\colon S\to \RR^+$. Для  произвольной ${\mathfrak m}_1$-интегрируемых  функций  $v$  на окружностях $\partial D_z\bigl(r(z)\bigr)$ при $z\in S$
можно определить 
{\it интегральные средние  с переменным радиусом $r$  по  окружностям}
\begin{equation}
v^{\circ r}\colon z\underset{z\in S}{\longmapsto}
\frac{1}{2\pi r(z)}\int_{\partial \overline  D_z(r(z))}v\dd {\mathfrak m}_1 
=
\frac{1}{2\pi} \int_{0}^{2\pi}  v\bigl(z+r(z)e^{i\theta}\bigr) \dd \theta \in \overline \RR.
\label{vpC}
\end{equation} 
Если функция $v$ определена на объединении кругов 
\begin{equation}\label{Sdrcup}
 S^{\cup r}:=\bigcup_{z\in S}\overline  D_z\bigl(r(z)\bigr)\subset \mathbb C, 
\end{equation}
то можем определить её {\it точную верхнюю грань по кругам }
\begin{equation}\label{veev}
v^{\vee r}\colon z\underset{z\in S}{\longmapsto} \sup_{\overline D_z(r(z))} v\in \overline \RR,
\end{equation}
а также {\it интегральные средние  с переменным радиусом $r$   по кругам\/}
\begin{equation}
v^{\bullet r}\colon z\underset{z\in S}{\longmapsto}
\frac{1}{\pi (r(z))^2} \int_{\overline  D_z(r(z))}  v \dd \mathfrak m_2  
\label{vpD}
\end{equation}
при  $\mathfrak m_2$-интегрируемости  $v$ на кругах $\overline  D_z\bigl(r(z)\bigr)$ при $z\in S$. 

Для любой функции $u$, субгармонической  на открытой окрестности объединения кругов 
$ S^{\cup r}$ из  \eqref{Sdrcup}, имеем неравенства \cite[теорема 2.6.8]{Rans}
\begin{equation}\label{vbc}
u\leq u^{\bullet r}\leq u^{\circ r} \leq u^{\vee r}
\quad\text{на  } S.
\end{equation}

\textit{Субгармоническую на $\CC$ функцию $u$  конечного типа\/} $\type[u]<+\infty$ (при порядке $1$)
с \textit{конечным логарифмическим интегралом\/} 
$J[u^+]\overset{\eqref{Jm0}}{<}+\infty$
называем \textit{субгармонической функцией класса Картрайт.\/}
В \cite[{\bf 3}, определение]{MS} и 
\cite[1.3.1]{BaiTalKha16}  так назывался существенно  \textit{более узкий класс\/} $\mathcal C$
субгармонических функций $u$ на $\CC$ конечного типа при порядке $1$, удовлетворяющих  условиям 
\textit{гармоничности на\/} $\CC\setminus \RR$, \textit{зеркальной симметричности относительно\/} $\RR$, т.е. $u(z)=u(\Bar z)$ для всех $z\in \CC$, а также 
\textit{с нулевым значением\/  $u(0)=0$ в нуле} и \textit{конечным интегралом}
$$
\int_{-\infty}^{+\infty} \frac{u^+(x)}{x^2}\dd x<+\infty.
$$

Целая функция $f$ экспоненциального  типа 
называется  \textit{целой  функцией класса Картрайт} \cite{Levin56}, \cite{Levin96}, \cite{HJ94}, 
если  $u:=\ln|f|$  --- субгармоническая функция класса Картрайт, т.е. $J\bigl[\ln^+|f|\bigr]\overset{\eqref{Jm0}}{<}+\infty$.

Наш основной результат этой первой части работы ---  следующая 
 
\begin{theorem}\label{th1} Если   $r\colon \CC\to (0,1]$ --- полунепрерывная снизу функция, для которой  
\begin{equation}\label{r}
\liminf\limits_{z\to \infty} \dfrac{\ln r(z)}{\ln |z|}>-\infty,
\end{equation} 
а  $u\neq -\infty$ и $M\neq -\infty$ ---  пара  субгармонических функций 
класса Картрайт, то  для любых  чисел $a>\type[u]$,  $c>\type[M]$,  и $d\in (0,2]$ существует  ограниченная на $\RR$  целая функция $h\neq 0$ экспоненциального типа  $\type_h\leq c$, для которой
\begin{equation}\label{{uvh}s}
\bigl(u(z)-M^{\bullet  r}(z)\bigr)+\ln \bigl|h(z)\bigr|\leq 
a|\Im z|\quad\text{при всех $z\in \CC$},
\end{equation}
где согласно  \eqref{vbc} функцию $M^{\bullet  r}$ можно заменить на 
$M^{\circ r} $ или $M^{\vee r}$ из \eqref{vpC}--\eqref{veev}.
При этом  найдётся такое исключительное множество $E_r\subset \CC$, что
 \begin{equation}\label{{uvh}E}
u(z)+\ln \bigl|h(z)\bigr|\leq M(z)+a |\Im z|
\quad\text{при всех $z\in \CC\setminus E_r$,}
\end{equation}
и в то же время   $d$-мерные обхваты Хаусдорфа множества $E_r$ 
переменного радиуса $r$ имеют ограничения 
\begin{equation}\label{{uvh}m}
{\mathfrak m}_d^{r}(E_r\cap S)\leq \sup_{z\in S} r(z) \quad\text{для  любого  $S\subset  \CC$}.
\end{equation}
В частности, при  $d:=1$ найдётся    исключительное  множество $Y_r\subset \RR$, для которого выполнено 
неравенство 
\begin{equation}\label{mE}
{\mathfrak m}_1\bigl(Y_r\cap (\RR \setminus [-y,y])\bigr)
\leq 2\sup\limits_{|\Im z|>y}r(z)\quad\text{ для  любого  $y\in  \RR^+$},
\end{equation}
а на прямых  $\bigl\{x+iy\bigm| x\in \RR \bigr\}$, параллельных $\RR$ и не проходящих через $iY_r$,   имеют место неравенства 
\begin{equation}\label{{uvh}Y}
u(x+iy)+\ln \bigl|h(x+iy)\bigr|\underset{x\in \RR}{\leq} M(x+iy)+
a |y|\quad\text{при  каждом   $y\in \RR\setminus Y_r$}.
\end{equation}
\end{theorem}

\begin{example}
При любом  $P\in\RR^+$ для непрерывной функции 
\begin{equation}\label{rP}
r\colon z\underset{z\in \CC}{\longmapsto}\dfrac{1}{(1+P+|z|)^P}\in (0,1]
\end{equation} 
имеем \eqref{r}, поскольку в случае  \eqref{rP} (нижний) предел из  \eqref{r}
равен $-P>-\infty$. 
\end{example}

\section{Доказательство основного результата}

Для функции $u\colon \RR\to \overline \RR$ 
при ${\mathfrak m}_1$-интегрируемости на $\RR$ функции 
\begin{equation}\label{u2}
u_{\pi}\colon x\underset{x\in \RR}{\longmapsto} \frac{1}{\pi}\dfrac{u(x)}{1+x^2}
\end{equation}
при всех $z\in \CC\setminus \RR$ через \textit{ядро Пуассона}
\begin{equation}\label{kP}
{\mathrm P}_{\CC\setminus \RR}\colon (z,x)\underset{z\in \CC\setminus \RR, x\in \RR}{\longmapsto} 
\frac{1}{\pi}\Bigl|\Im \frac{1}{z-x}\Bigr|\in \RR^+
\end{equation} 
на  $\CC\setminus \RR$ определён \textit{интеграл Пуассона\/} 

\begin{equation}\label{Poi}
{\mathcal  P}_{\CC\setminus \RR} u
\colon z\underset{z\in \CC\setminus \RR}{\longmapsto}
\int_{\RR} {\mathrm P}_{\CC\setminus \RR} (z,x)u(x)\dd x
=
\frac{1}{\pi}\int_{-\infty}^{+\infty} 
\frac{|\Im z| u(x)}{(\Im z)^2+(\Re z-x)^2}\dd x\in \overline \RR.
\end{equation}
Допустим, что  функция $u_{\pi}$, определённая в \eqref{u2},   ${\mathfrak m}_1$-{\it суммируема,\/} 
т.е. конечен интеграл $J\bigl[|u|\bigr]\in \RR$. Тогда интеграл  Пуассона \eqref{Poi} определяет  гармоническую функцию на $\CC\setminus \RR$, которая называется  \textit{гармоническим продолжением функции\/ $u$ на\/ $\CC\setminus \RR$, {\rm т.е.} вне\/ $\RR$,\/}  
а если   функция $u$ при этом \textit{ непрерывна на\/} $\RR$, то функция 
\begin{equation}\label{limP}
z\underset{z\in \CC}{\longmapsto} 
\begin{cases}
\bigl({\mathcal  P}_{\CC\setminus \RR} u\bigr)(z)&\text{при $z\in \CC\setminus \RR$},\\
 u(z) &\text{при $z\in \RR$}
\end{cases}
\end{equation}
\textit{непрерывна на\/  $\CC$.}
При  $J\bigl[|u|\bigr]\in \RR$ и более слабом условии лишь \textit{полунепрерывности сверху\/}  функции 
$u$ на $\RR$,  
функция \eqref{limP}  \textit{полунепрерывная сверху на\/  $\CC$,} поскольку функцию $u$ на $\RR$ нетрудно представить как предел убывающей последовательности непрерывных функций  $u_n$ с  
${\mathfrak m}_1$-суммируемыми  на $\RR$ функциями  $(u_n)_{\pi}$, определёнными  в \eqref{u2}.  
При этом из вида \eqref{Jm0} логарифмического интеграла,  \eqref{kP} ядра Пуассона  и по определению \eqref{Poi} интеграла Пуассона имеем
\begin{equation}\label{pmi}
J[u]=\bigl({\mathcal  P}_{\CC\setminus \RR} u\bigr)(i)=\bigl({\mathcal  P}_{\CC\setminus \RR} u\bigr)(-i).  
\end{equation}

{\it Сдвиг   функции $v$ на $z_0\in \CC$} и  её \textit{гомотетия с коэффициентом}  $k\in \RR$ --- 
это соответственно  функции $v(\cdot-z_0)\colon  z\longmapsto v(z-z_0)$
и $v(k\cdot)\colon  z\longmapsto v(kz)$.  

\begin{lemma}\label{lem1} Если  $u\neq -\infty$ --- субгармоническая функция класса Картрайт,
то конечен логарифмический интеграл $J\bigl[(-u)^+\bigr] \overset{\eqref{Jm0}}{<}+\infty$ 
от положительной части $(-u)^+$ противоположной функции $-u$,\/
 ${\mathfrak m}_1$-суммируема функция  $u_{\pi}$ из \eqref{u2} и 
\begin{subequations}\label{IP}
\begin{align}
u(z)&\underset{z\in \CC}{\leq} u^{\bal}(z)\underset{z\in \CC}{:=}
\bigl({\mathcal  P}_{\CC\setminus \RR} u\bigr)(z) +\type[u]|\Im z|,
\tag{\ref{IP}b}\label{{IP}b}
\\
J\bigl[u(\cdot-iy_0)\bigr]&\underset{y_0\in \RR}{\leq} 
\bigl({\mathcal  P}_{\CC\setminus \RR} u\bigr)(i+i|y_0|\bigr)
+\type[u]|y_0|,
\tag{\ref{IP}J}\label{{IP}J}
\\
J\bigl[u^+(\cdot-iy_0)\bigr]&\underset{y_0\in \RR}{\leq} 
\bigl({\mathcal  P}_{\CC\setminus \RR} u^+\bigr)(i+i|y_0|\bigr)
+\type[u]|y_0|,
\tag{\ref{IP}J$^+$}\label{{IP}J+}
\end{align}
\end{subequations}
где $u^{\bal}\neq -\infty$ из \eqref{{IP}b}  --- субгармоническая функция класса Картрайт, 
гармоническая на $\CC\setminus \RR$, с совпадающими   сужениями  $u^{\bal}{\lfloor}_\RR=u{\lfloor}_\RR$  на $\RR$.

В частности, сдвиг субгармонической функции класса Картрайт на любое $z_0\in \CC$ 
--- субгармонической функции класса Картрайт, а гомотетия субгармонической функции класса Картрайт с любым коэффициентом $k\in \RR$ также даёт субгармоническую функцию класса Картрайт, или, более детально, 
\begin{equation}\label{ku}
J\bigl[u(k\cdot)\bigr]=\bigl({\mathcal  P}_{\CC\setminus \RR} u\bigr)(\pm ik\bigr), \quad J\bigl[u^+(k\cdot)\bigr]=\bigl({\mathcal  P}_{\CC\setminus \RR} u^+\bigr)(\pm ik\bigr).
\end{equation}
\end{lemma}
\begin{proof} 
Для  субгармонической функции $u\not\equiv -\infty$ действие  на неё   {\it оператора Лапласа\/} ${\bigtriangleup}$  в смысле теории обобщённых функций  определяет её {\it распределение масс Рисса\/}
$\frac{1}{2\pi}{\bigtriangleup}u=:\varDelta_u$.
Для любой субгармонической функции $u$ конечного типа согласно 
\cite[предложение 4.1, (4.19)]{KhaShmAbd20}, \cite[лемма 2]{SalKha20U}
существует число $C\in \RR^+$, с  которым   имеют место неравенства 
\begin{multline*}
\max \biggl\{
\int_{\CC^{\up}\cap ((R\overline \DD)\setminus \overline \DD)} \Big|\Im \frac{1}{z}\Big| \dd \varDelta_u(z), 
\int_{(-\CC^{\up})\cap ((R\overline \DD)\setminus \overline \DD)} \Big|\Im \frac{1}{z}\Big| \dd \varDelta_u(z)
\biggr\}
\\
\leq \frac{1}{2\pi}\int_1^R\frac{u(x)+u(-x)}{x^2}\dd x +C
\quad\text{при всех $R>1$,}
\end{multline*}
где левая часть положительна. Отсюда по представлению  $u:=u^++(-u)^+$ 
\begin{multline*}
\frac{1}{2\pi}\int_1^R\frac{(-u)^+(x)+(-u)^+(-x)}{x^2}\dd x\leq  \frac{1}{2\pi}\int_1^R\frac{u^+(x)+u^+(-x)}{x^2}\dd x +C\\
\leq \frac{1}{\pi}\int_1^R\frac{u^+(x)+u^+(-x)}{1+x^2}\dd x
+C\leq 2J[u^+]+C \quad\text{при всех $R>1$,}
\end{multline*}
после чего, устремляя в левой части  $R$ к $+\infty$, получаем $J\bigl[(-u)^+\bigr]<+\infty$ и  
$J\bigl[|u|\bigr]=J[u^+]+J\bigl[(-u)^+\bigr]<+\infty$, т.е. 
${\mathfrak m}_1$-суммируема функция  $u_{\pi}$ из \eqref{u2}.

В \cite[1.2.2, \S~6]{KhaShm19} для субгармоническая функции $u$ конечного типа при $J[u^+]<+\infty$ конструируется \textit{субгармоническое выметание\/}  
одновременно из верхней полуплоскости $\CC^{\up}$ и нижней полуплоскости $-\CC^{\up}$ как субгармоническая функция $u^{\bal}$ из  \eqref{{IP}b}, гармоническая на $\CC\setminus \RR$, а также с совпадающими   сужениями  $u^{\bal}{\lfloor}_\RR=u{\lfloor}_\RR$  на $\RR$. 
Далее  из  \eqref{{IP}b} при всех $z\in \CC$ и $y_0\in \RR$ получаем
\begin{equation}\label{Pxy+}
u(z-iy_0)\overset{\eqref{{IP}b}}{\leq} (u)^{\bal} (z-iy_0)
=\bigl({\mathcal  P}_{\CC\setminus \RR} u\bigr)(z-iy_0) 
+\type[u]|\Im z-y_0|,
\end{equation}

 Слева в \eqref{Pxy+} стоит значение в точке $z$ субгармонической функции 
$u(\cdot-iy_0)$, полученной как сдвиг в точку $iy$ субгармонической  функции $u$.
Применение  \eqref{pmi} к сдвигу $u(\cdot-iy_0)$ даёт соотношения
\begin{multline}\label{ineqa}
J\bigl[u(\cdot -iy_0)\bigr]\overset{\eqref{pmi}}{=}
\bigl({\mathcal  P}_{\CC\setminus \RR} u(\cdot-iy_0)\bigr)(\pm i)
\\ 
\overset{\eqref{Pxy+}}{\leq}  \biggl({\mathcal  P}_{\CC\setminus \RR} 
\Bigl({\mathcal  P}_{\CC\setminus \RR} \bigl(u(\cdot-iy_0)\bigr) +
\type[u]|\Im \cdot -y_0|\Bigr)
\biggr)(\pm i)
\\
\overset{\eqref{pmi}}{=}
\biggl({\mathcal  P}_{\CC\setminus \RR}\Bigl({\mathcal  P}_{\CC\setminus \RR} \bigl(u(\cdot-iy_0)\bigr) \Bigr)\biggr)(\pm i)
+\type[u]|y_0|.
\end{multline}
Для первого слагаемого в правой части  \eqref{ineqa} ввиду гармоничности функции 
$\bigl({\mathcal  P}_{\CC\setminus \RR} u\bigr)(\cdot-iy_0)$ в нижней полуплоскости $-\CC^{\up}$
имеют место равенства 
$$
\biggl({\mathcal  P}_{\CC\setminus \RR}\Bigl({\mathcal  P}_{\CC\setminus \RR} \bigl(u(\cdot-iy_0)\bigr) \Bigr)\biggr)(- i)
=\Bigl(\bigl({\mathcal  P}_{\CC\setminus \RR} u\bigr)(\cdot-iy_0) \Bigr)(-i) 
=
\bigl({\mathcal  P}_{\CC\setminus \RR} u\bigr)(-i-iy_0) , 
$$
Последнее, подставленное в правую часть  \eqref{ineqa},
даёт 
$$
J\bigl[u(\cdot -iy_0)\bigr]\overset{\eqref{ineqa}}{\leq} 
\bigl({\mathcal  P}_{\CC\setminus \RR} u\bigr)(-i-iy_0) 
 +\type[u]|y_0|,
$$
откуда ввиду зеркальной симметрии значений интеграла Пуассона \eqref{Poi} относительно вещественной оси 
получаем требуемое \eqref{{IP}J}.
Функция $u^+$ также \textit{субгармоническая\/} функция, и очевидно, класса Картрайт с $\type[u^+]
\overset{\eqref{typevf}}{=}\type[u]$. Поэтому \eqref{{IP}J} влечёт за собой \eqref{{IP}J+}. 
В частности,  сдвиг $u(\cdot-iy_0)$ при любом $y_0\in \RR$ --- субгармоническая функция класса Картрайт. 
Сдвиг $u(\cdot -x_0)$ на $x_0\in \RR$ тоже субгармоническая функция класса Картайт, 
поскольку 
\begin{equation}\label{Jux0}
J\bigl[u^+(\cdot-x)\bigr]\overset{\eqref{Jm0}}{=}\frac{1}{\pi}
 \int_{-\infty}^{+\infty}\frac{u^+(x-x_0)}{1+x^2}\dd x
\leq J[u^+]\sup_{x\in \RR}\frac{1+(x-x_0)^2}{1+x^2}<+\infty.
\end{equation}
Таким образом, сдвиг $u(\cdot -z_0)$ на $z_0=x_0+iy_0\in \CC$ --- субгармоническая функция класса Картрайт как результат последовательного применения двух сдвигов на $iy_0$ и $x_0$.
Наконец, при гомотетии с  коэффициентом $k\neq 0$, используя замену переменной, получаем 
\begin{equation*}
J\bigl[u(k\cdot)\bigr]\overset{\eqref{Jm0}}{=}\frac{1}{\pi}
 \int_{-\infty}^{+\infty}\frac{u(kx)}{1+x^2}\dd x
=\frac{1}{\pi}
\int_{-\infty}^{+\infty}\frac{|k|u(t)}{k^2+t^2}\dd t
\overset{\eqref{Poi}}{=}\bigl({\mathcal  P}_{\CC\setminus \RR} u\bigr)(\pm ik),
\end{equation*}
а для  гомотетии с коэффициентом $0$ равенство \eqref{ku} тривиально.
\end{proof}

{\it Горизонтальную  открытую полосу 
ширины\/ $2b$, симметричную  относительно вещественной оси\/ $\RR$,} обозначаем через 
$\strip_b:=\Bigl\{z\in \CC\Bigm| |\Im z|< b\Bigr\}$.

\begin{lemma}\label{lem2} При любых $b\in \RR^+\setminus 0$ и $y\in \RR$ для любой  субгармонической функции $u\neq -\infty$  класса Картрайт найдётся целая функция $F\neq 0$ класса Картрайт типа  $\type_F=\type[u]$, 
для которой $\ln\bigl|F\bigr|\geq u$ вне полосы   $iy+\strip_b$, т.е.
\begin{equation}\label{outst}
u(z)\leq \ln\bigl| F(z)\bigr|\quad\text{при $|\Im z-y|\geq b$.}
\end{equation}
\end{lemma}
\begin{proof} По лемме  \ref{lem1} сдвиг функции не выводит её из класса Картрайт, поэтому достаточно рассмотреть случай $y=0$. Пусть $u^{\bal}\geq u$ --- субгармоническая функция класса Картрайт из 
\eqref{{IP}b}. Достаточно построить  целую функцию $F\neq 0$ класса Картрайт  
для $u^{\bal}$, для которой выполнено \eqref{outst} с $u^{\bal}$ вместо $u$. 
Вследствие этого можно изначально считать, что \textit{ исходная субгармоническая функция  $u$
класса Картрайт гармоническая на $\CC\setminus \RR$. }

Теперь воспользуемся очень частным случаем одного результата, вытекающего из KKK
(Kjellberg-Kennedy-Katifi) аппроксимационного метода 
\cite[10.5]{Hayman2}. 

\begin{lemma}[{\!\!\rm\cite[Лемма 10.12]{Hayman2}, \cite[лемма 2.1]{Kha01l}}]\label{lem:1} 
Пусть $u$ --- субгармоническая фу\-н\-к\-ц\-ия конечного типа, гармоническая на $\CC\setminus \RR$. Тогда существует такая целая функция $g$ экспоненциального типа с нулями только
на $\RR$, что
\begin{equation}\label{gu}
\ln \bigl|g(z)\bigr| = u(z) +O\Bigl(\ln^+\frac1{|\Im z|}\Bigr)
+O(\ln |z|) \quad\text{при  $z\to \infty$}.
\end{equation}

\end{lemma}
Для целой функции $g\neq 0$ экспоненциального типа из леммы \ref{lem:1}
по  соотношению \eqref{gu}  для любого $b\in \RR^+\setminus 0$
 существует $C\in \RR^+$, для которого 
\begin{equation}\label{pgu}
\Bigl|u(z)-\ln\bigl|g(z)\bigr|\Bigr|\overset{\eqref{gu}}{\leq} C\ln\bigl(2+|z|\bigr)
\quad\text{при всех $z\in \CC\setminus \strip_b$,}
\end{equation}   
а также $\type_g=\type[u]$. В частности, из неравенств 
$$
\ln|g(x-ib)|\overset{\eqref{pgu}}{\leq} u(x-ib)+ C\ln \bigl(1+\sqrt{x^2+b^2}\bigr)
\quad\text{при всех $x\in \RR$,}
$$
 по лемме   \ref{lem1} следует, что $g$ --- целая функция класса Картрайт. При этом из \eqref{pgu}
для достаточно большого $N\in \NN$ получаем 
$u(z)\overset{\eqref{gu}}{\leq} \ln\bigl|Nz^Ng(z)\bigr|$
при всех $z\in \CC\setminus \strip_b$,
откуда  целая функция $F\colon z\underset{z\in \CC}{\longmapsto}Nz^Ng(z)$
 и есть требуемая.
\end{proof}

\begin{lemma}[{очень частный случай \cite[основная теорема]{Kha22arxiv}, \cite[теорема 9]{Kha22IzRAN}}]\label{M}
Пусть $M\neq -\infty$ --- субгармоническая фу\-н\-к\-ция конечного типа на $\CC$, полунепрерывная снизу  функция 
 $r\colon \CC\to (0,1]$ удовлетворяет условию \eqref{r}, а также $d\in (0,2]$.
 Тогда существуют целая функция $f\neq 0$ экспоненциального типа $\type_f\leq \type[M]$ и исключительное множество $E_r\subset \CC$, удовлетворяющее \eqref{{uvh}m}, для которых 
\begin{equation}\label{ubc}
\ln |f|\overset{\eqref{vpD}}{\leq} M^{\bullet r} 
\quad\text{на  $\CC$},\qquad \ln |f|\leq M \quad\text{на  $\CC\setminus E_r$.}
\end{equation}
\end{lemma}

\begin{lemma}\label{lem5} Если в лемме\/ {\rm \ref{M}}   функция $M$ класса Картрайт и $d:=1$, то  
 целая функция $f\neq 0$ из \eqref{ubc} класса Картрайт  типа $\type_f\leq \type[M]$ и 
найдётся такое   исключительное множество $Y_r\subset \RR$,
удовлетворяющее  \eqref{mE}, что 
\begin{equation}\label{lmpr}
\ln \bigl|f(x+iy)\bigr|
\underset{x\in \RR}{\leq} M(x+iy)\quad\text{при  каждом   $y\in \RR\setminus Y_r$}.
\end{equation}
\end{lemma}
\begin{proof}  При $d:=1$ по определению \eqref{mr} 
одномерного обхвата Хаусдорфа с радиусом обхвата $r$ 
ортогональная проекция $iY_r\subset i\RR$  на мнимую ось $i\RR$
множества $E_r$, удовлетворяющего \eqref{{uvh}m}, как нетрудно видеть,  
удовлетворяет соотношениям  \eqref{mE} и при каждом $y\in \RR\setminus Y_r$ прямые $\bigl\{x+iy\bigm| x\in \RR\bigr\}$ не пересекают множество  $E_r$. 
Тогда по  второму неравенству в \eqref{ubc} получаем \eqref{lmpr}, откуда  по  лемме  \ref{lem1}
целая функция $f$ является  целой функцией класса Картрайт.
\end{proof}

\begin{proof}[Доказательство теоремы\/ {\rm \ref{th1}}]
По лемме \ref{lem2} найдётся  целая функция $F\neq 0$  класса Картрайт и типа $\type_F\leq \type[u]$, для которой  выполнено \eqref{outst} при некотором выборе $0<b<y\in \RR^+$. 
Рассмотрим  строго положительное число 
\begin{equation}\label{bq}
 q:=\frac12\min \bigl\{a-\type[u], c-\type[M]\bigr\}>0. 
\end{equation}
По теореме Бёрлинга\,--\,Мальявена существует \textit{ограниченная на $\RR$
целая функция\/ $h_F\neq 0 $   экспоненциального типа}
\begin{equation}\label{thF}
\type_{h_F}\leq q \overset{\eqref{bq}}{\leq}\frac{1}{2}\bigl(a-\type[u]\bigr),
\end{equation} 
для которой целая функция $Fh_F$   экспоненциального  типа 
\begin{equation}\label{thFcd}
\type_{Fh_F}\overset{\eqref{thF}}{\leq} \type[u]+\frac{1}{2}\bigl(a-\type[u]\bigr)
\end{equation} удовлетворяет неравенствам 
\begin{equation}\label{uhF}
\ln|F(x)|+\ln \bigl|h_F(x)\bigr|\overset{\eqref{uh}}{\leq} 0 \quad\text{при всех $x\in \RR$}.  
\end{equation}
Отсюда для субгармонической функции $\ln|Fh_F|$ по лемме \ref{lem1} получаем 
\begin{multline*}
\bigl(\ln|Fh_F|\bigr)(z)\overset{\eqref{{IP}b}}{\underset{z\in \CC}{\leq}}
\bigl({\mathcal  P}_{\CC\setminus \RR} \ln|Fh_F|\bigr)(z) +\type\bigl[\ln|Fh_F|\bigr]|\Im z|\\
\overset{\eqref{uhF},\eqref{thF}}{\leq} 
\Bigl(\type[u]+\frac{1}{2}\bigl(a-\type[u]\bigr)\Bigr)|\Im z|,
\end{multline*}
что согласно \eqref{outst} при  $|\Im z-y|\geq b$ даёт неравенство 
\begin{equation}\label{aqqq}
u(z)+\ln\bigl|h_F(z)\bigr|\leq 
\Bigl(\type[u]+\frac{1}{2}\bigl(a-\type[u]\bigr)\Bigr)|\Im z|.
\end{equation}
В правой части  \eqref{aqqq} функция гармоническая на  $\CC\setminus \RR$ и, исходя из выбора  чисел $0<b<y$, 
неравенство \eqref{aqqq}, справедливое всюду  вне полосы $iy+\strip_b\subset \CC^{\up}$, можно, используя 
субгармонический  вариант теоремы Фрагмена\,--\,Линделёфа для полос, продолжить на все точки $z\in \CC$. 

По лемме \ref{lem5} в сочетании с леммой \ref{M} существуют целая функция $f\neq 0$ класса Картрайт типа $\type_f\leq \type[M]$ и исключительное множество $E_r\subset \CC$, удовлетворяющее \eqref{{uvh}m}, для которых 
выполнено \eqref{ubc},  а при выборе  $d:=1$ для некоторого исключительного множества $Y_r\subset \RR$
имеем ещё и \eqref{mE} вместе с  \eqref{lmpr}. 
По теореме Бёрлинга\,--\,Мальявена существует ограниченная на $\RR$
целая функция $h_f\neq 0 $ экспоненциального  типа
\begin{equation}\label{thFM}
\type_{h_f}\leq q \overset{\eqref{bq}}{=} 
\frac12 \min \bigl\{c-\type[M], a-\type[u]\bigr\},
\end{equation} 
для которой целая функция $fh_f$ экспоненциального типа 
\begin{equation}\label{thFcdM}
\type_{fh_f}\overset{\eqref{efet}}{\leq} \type_f+\type_{h_f}
\overset{\eqref{thFM}}{\leq} \type[M]+q
\end{equation} 
удовлетворяет неравенствам 
\begin{equation}\label{uhFM}
\ln|f(x)|+\ln \bigl|h_f(x)\bigr|\overset{\eqref{uh}}{\leq} 0 \quad\text{при всех $x\in \RR$}.  
\end{equation}
При этом домножая, при необходимости, функцию $h_f$ на достаточно малое строго положительное число и сохраняя за произведением прежнее обозначение $h_f$, можем добиться с сохранением \eqref{uhFM} того, что 
$|h_f|\leq 1$ на $\RR$, откуда    по лемме \ref{lem1}, применённой 
к субгармонической функции $\ln|h_f|$, получаем 
\begin{equation}\label{h0}
\ln \bigl|h_f(z)\bigr|\overset{\eqref{{IP}b}}{\leq} 
\bigl({\mathcal  P}_{\CC\setminus \RR} \ln|h_f|\bigr)(z) +\type\bigl[\ln|h_f|\bigr]|\Im z|
\overset{\eqref{thFM}}{\leq} 
q|\Im z|.
\end{equation}
Далее, складывая неравенство  \eqref{aqqq} с первым неравенством в \eqref{ubc} и с 
крайними частями неравенств \eqref{h0},
получаем при всех $z\in \CC$ неравенства 
\begin{multline}\label{aMM}
u(z)+\ln\bigl|h_F(z)\bigr|+\ln\bigl|f(z)\bigr|+
\ln\bigl|h_f(z)\bigr|
\\
\overset{\eqref{aqqq},\eqref{ubc},\eqref{h0}}{\leq} 
\Bigl(\type[u]+\frac{1}{2}\bigl(a-\type[u]\bigr)\Bigr)|\Im z|
+M^{\bullet r}(z)+q|\Im z|\\
\leq M^{\bullet r}(z) +
\Bigl(\type[u]+\frac{1}{2}\bigl(a-\type[u]\bigr)+q\Bigr)|\Im z|
\overset{\eqref{thFM}}{\leq}
M^{\bullet r}(z) +a|\Im z|,
\end{multline}
что  для целой функции  
\begin{equation}\label{hhh}
h:=h_Ffh_f
\end{equation} 
можно записать как неравенство \eqref{{uvh}s} при всех $z\in \CC$. 
При этом функция $h$ ограничена на $\RR$ как произведение \eqref{hhh}
ограниченной  на $\RR$ согласно \eqref{uhFM} функции $fh_f$
на  ограниченную на $\RR$ функцию $h_F$, что отмечено выше перед \eqref{thF}.
Наконец, оценка типа целой функции $h$ следует из неравенств 
\begin{multline*}
\type_h\overset{\eqref{hhh}}{=}\type_{h_Ffh_f}
\overset{\eqref{efet}}{\leq} 
\type_{h_F}+\type_{fh_f}
\overset{\eqref{thF},\eqref{thFcdM}}{\leq}
q+\bigl(\type[M]+q\bigr)\\
=\type[M]+2q\overset{\eqref{bq},\eqref{thFM}}{\leq} 
\type[M]+2\cdot \frac12\bigl(c-\type[M]\bigr)=c.
\end{multline*}
Таким образом,  функция $h$ для \eqref{{uvh}s} с требуемыми свойствами построена. 
При этом по лемме \ref{M} из  второго неравенства  
для той же функции $h$ имеем \eqref{{uvh}E} вкупе с \eqref{{uvh}m}. 
Окончательно в случае $d:=1$ лемма \ref{lem5} с неравенствами \eqref{lmpr} 
обеспечивает для $h$ выполнение \eqref{mE}--\eqref{{uvh}Y}, и теорема \ref{th1} доказана.  
\end{proof}

\begin{rem}
В продолжение настоящей статьи  намечается более детально исследовать субгармонические функции класса Квртрайт, 
на основе чего будут   рассмотрены субгармонические аналоги теоремы Бёрлинга\,--\,Мальявена 
о радиусе полноты \cite{BM67}, \cite{Kah62}, \cite{Red77}, 
\cite{Kra89}, \cite{Koosis92}, \cite{Kha94}, \cite{Koosis96},  \cite{HJ94}, \cite{Khsur}, 
\cite{KhaTalKha14}, \cite{BaiTalKha16}.
\end{rem}

\bigskip

Адреса авторов:

   \medskip

\noindent
    Б.Н. ХАБИБУЛЛИН

\noindent
Институт математики с вычислительным центром УФИЦ РАН,
     
\noindent
      450008, г. Уфа, ул. Чернышевского, 112,

\noindent
        khabib-bulat@mail.ru

\medskip

\noindent
    Е.Г. КУДАШЕВА

\noindent
Башкирский государственный педагогический университет им. М. Акмуллы,
     
\noindent
      450008, г. Уфа, ул. Октябрьской революции, 3А,

\noindent
        lena\_kudasheva@mail.ru


\bigskip
\begin{thebibliography}{86}%% В порядке цитирования

\bibitem{Boas}
 R.\,P.~Boas, Jr.,  {\it Entire Functions,\/}   New York:  Academic Press,   1954.  

\bibitem{Levin96} 
B.~Ya.~ Levin, {\it Lectures on entire functions,\/} 
Transl. Math. Mono\-graphs, {\bf 150}, Providence RI:  Amer. Math. Soc.,  1996.

\bibitem{RC}
 L.\,A.~Rubel (with J.\,E. Colliander) {\it Entire and Meromorphic Functions,\/}
  New York--Berlin--Heidelberg:  Sprin\-ger-Verlag, 1996.


\bibitem{Khsur}
 Б.~Н.~Хабибуллин, {\it Полнота систем экспонент и множества единственности,\/}  
издание  четвёртое дополненное, {\bf 2}, Уфа:  РИЦ БашГУ, 2012.
%%\finalinfo 
  \href{https://www.researchgate.net/publication/271841461}{https://www.researchgate.net/publication/271841461}


\bibitem{Levin56}
 Б.\,Я.~Левин, {\it Распределение корней целых функций,\/}
  М.: ГИТТЛ, 1956.

\bibitem{HJ94}
 V.~Havin, B.~J\"oricke, {\it The Uncertainty Principle in Harmonic Analysis,\/}
  Berlin: Springer-Verlag, 1994.

\bibitem{Koosis88}
 P.~Koosis, {\it The logarithmic integral.\/} I, Cambridge Stud. Adv. Math.,  {\bf 12},  
  Cambridge:  Cambridge Univ. Press, 1988.

\bibitem{Koosis92}
 P.~Koosis,   {\it The logarithmic integral.\/} II,
Cambridge Stud. Adv. Math.,  {\bf 21},  Cambridge:  Cambridge Univ. Press, 1992.

\bibitem{Koosis96}
 P.~Koosis,    {\it Le\c{c}ons sur le th\'eor\`eme de Beurling et Malliavin,\/}
Univ. Montr\'eal, Montr\'eal, QC: Les Publications CRM, 1996. 


\bibitem{MasNazHav05}
 Дж.~Машреги, Ф.~Л.~Назаров, В.~П.~Хавин,
 {\it Теорема Бёрлинга\,--\,Мальявена о~мультипликаторе: седьмое доказательство,\/}
 Алгебра и анализ, {\bf 17}:5 (2005), 3--68.

\bibitem{BM62} 
 A.~Beurling A., P.~Malliavin,
{\it  On Fourier transforms of measures with compact support,\/}
  Acta Math., {\bf 107} (1962),  1291--309.

\bibitem{Mal79} 
 P.~Malliavin, {\it On the multiplier theorem for Fourier transforms of measures with compact support,\/}
   Ark. Mat., {\bf 17}  (1979). 69--81. 

\bibitem{Kiselman}
 Ch.\,O.~Kiselman, {\it Order and type as measures of growth for convex or entire functions,\/}
  Proc.  London Math. Soc., {\bf 66}:3 (1993), 152--86.


\bibitem{KhaShm19}
 Б.~Н.~Хабибуллин, А.~В.~Шмелёва, 
{\it  Выметание мер и субгармонических функций на систему лучей.\/ {\rm I.} Классический случай,\/}
  Алгебра и анализ, {\bf 31}:1  (2019), 156--210.


\bibitem{Kha22}
  Б.~Н.~Хабибуллин,
{\it Интегралы от разности субгармонических функций по мерам и характеристика Неванлинны,\/}
Математический  сборник, {\bf 213}:5 (2022), 126--166.

\bibitem{EG}
  Л. К. Эванс,  К. Ф. Гариепи, {\it Теория меры  и тонкие свойства функции,\/}
  Новосибирск:  Научная книга (ИДМИ), 2002.

\bibitem{Rans}
 Th.~Ransford, {\it Potential Theory in the Complex Plane,\/}
  Cambridge:  Cambridge University Press, 1995.

\bibitem{MS}
V. Matsaev, M. Sodin, 
{\it Distribution of Hilbert transforms of measures,\/}
Geom. Funct. Anal. {\bf 10}:1 (2000), 160--184.

\bibitem{BaiTalKha16}
 Т.~Ю.~Байгускаров, Г.~Р.~Талипова, Б.~Н.~Хабибуллин,
{\it Подпоследовательности нулей для классов целых функций экспоненциального типа, выделяемых ограничениями на их рост,\/}  Алгебра и анализ, 
{\bf 28}:2 (2016), 1--33.
%%\transl
%%  St. Petersburg Math. J.
%%  2017
%% 28
%% 2
%% 127--151

\bibitem{KhaShmAbd20}
 Б.~Н.~Хабибуллин, А.~В.~Шмелёва, З.~Ф.~Абдуллина,
 {\it Выметание мер и субгармонических функций на систему лучей.\/ {\rm II.} Выметания конечного рода и регулярность роста на одном луче,\/}  Алгебра и анализ,  {\bf 32}:1  (2020), 208--243.
%%\transl
%%  St. Petersburg Math. J.
%%  2021
%% 32
%% 1
%% 155--181


\bibitem{SalKha20U}
 А.~Е.~Салимова, Б.~Н.~Хабибуллин,
{\it Рост субгармонических функций вдоль прямой и распределение их  распределений масс  Рисса,\/}
  Уфимский математический журнал, {\bf 12}:2  (2020), 35--48.
%%\transl
%%  Ufa Math. J.
%%  2020
%% 12
%% 2
%% 35--49

\bibitem{Hayman2}
Hayman W.~K.,   {\it Subharmonic functions.\/} II, London:
Academic Press, 1989.

\bibitem{Kha01l}
 Б.~Н.~Хабибуллин, 
{\it О~росте целых функций экспоненциального типа с~нулями вблизи прямой,\/}
  Математические  заметки, {\bf 70}:4  (2001), 621--635.
%%\transl
%%  Math. Notes
%%  2001
%%% 70
%% 4
%% 560--573

\bibitem{Kha22arxiv}
B. N. Khabibullin, {\it The restriction from below of the subharmonic function by the logarithm of the module of entire function,\/} arXiv:2203.12383, 22 March 2022,  11 pages, in Russian  \href{https://arxiv.org/abs/2203.12383}{https://arxiv.org/abs/2203.12383}


\bibitem{Kha22IzRAN}
Б.Н. Хабибуллин, 
{\it Распределения корней и масс   целых и субгармонических функций с ограничениями на их рост  вдоль полосы,\/}
2022, 61 стр. (направлено в печать).

\bibitem{BM67} 
A.~Beurling A., P.~Malliavin, {\it On the closure of characters and the zeros of entire functions,\/}
 Acta Math., {\bf 118} (1967), 79--93.

\bibitem{Kah62}
J.-P.~Kahane, {\it Travaux de Beurling et Malliavin,\/}
S\'eminaire Bourbaki {\rm (ann\'ee 1961/62,  expos\'es 223--240, Talk no. 225)},   
no. 7  (1962), 27--39.

\bibitem{Red77} 
 R.~M.~Redheffer, {\it Completeness of sets of complex exponentials,\/}
Adv. in Math., {\bf 24} (1977), 1--62.

\bibitem{Kra89}
И.~Ф.~Красичков-Терновский, 
{\it Интерпретация теоремы Бёрлинга\,--\,Мальявена о~радиусе полноты,\/}
 Математический сборник, {\bf 180}:3 (1989), 397--423.

\bibitem{Kha94}
Б. Н. Хабибуллин, 
{\it Неконструктивные доказательства теоремы Бёрлинга\,--\,Мальявена о радиусе полноты и теоремы неединственности для целых функций,\/} Известия  РАН. Серия математическая, {\bf 58}:4 (1994), 125--148.

\bibitem{KhaTalKha14}
 Б. Н. Хабибуллин, Г. Р. Талипова, Ф. Б. Хабибуллин, 
{\it Подпоследовательности нулей для пространств Бернштейна и полнота систем экспонент в пространствах функций на интервале,\/} 
Алгебра и анализ, {\bf 26}:2 (2014), 185--215.

\end{thebibliography}
\end{document}